\documentclass[11pt]{article}
\usepackage{epic,latexsym,amssymb}
\usepackage{color}
\usepackage{tikz}

\newtheorem{theorem}{Theorem}
\newtheorem{lemma}{Lemma}

\newtheorem{cor}{Corollary}
\newtheorem{obser}{Observation}
\newtheorem{prop}{Proposition}
\newtheorem{conj}{Conjecture}

\newcommand{\QED}{$\Box$}

\newcommand{\modo}{{\rm mod}}
\newcommand{\cG}{\mathcal{G}}
\newcommand{\barD}{\overline{D}}
\newcommand{\pn}{{\rm pn}}
\newcommand{\IR}{{\rm IR}}

\newcommand{\coro}{{\rm cor}}
\newcommand{\epn}{{\rm epn}}

\newcommand{\dstart}{\gamma_{{\rm g}}}
\newcommand{\sstart}{\gamma_{{\rm g}}^\prime}

\newcommand{\cF}{{\cal F}}

\newcommand{\proof}{\noindent\textbf{Proof. }}
\newcommand{\2}{ \vspace{0.2cm} }
\newcommand{\1}{ \vspace{0.1cm} }

\let\oldenumerate\enumerate
\renewcommand{\enumerate}{
  \oldenumerate
  \setlength{\itemsep}{0pt}
  \setlength{\parskip}{0pt}
  \setlength{\parsep}{0pt}
}

\begin{document}

\title{The Enclaveless Competition Game}

\author{$^1$Michael A. Henning and $^2$Douglas F. Rall
\\ \\
$^1$Department of Mathematics and Applied Mathematics \\
University of Johannesburg \\
Auckland Park, 2006 South Africa\\
\small \tt Email: mahenning@uj.ac.za  \\
\\
$^2$Department of Mathematics \\
Furman University \\
Greenville, SC, USA\\
\small \tt Email: doug.rall@furman.edu}

\date{}
\maketitle

\begin{abstract}
For a subset $S$ of vertices in a graph $G$, a vertex $v \in S$ is an enclave of $S$ if $v$ and all of its neighbors are in $S$, where a neighbor of $v$ is a vertex adjacent to $v$. A set $S$ is enclaveless if it does not contain any enclaves. The enclaveless number $\Psi(G)$ of $G$ is the maximum cardinality of an enclaveless set in $G$. As first observed in 1997 by Slater [J. Res. Nat. Bur. Standards 82 (1977), 197--202], if $G$ is a graph with $n$ vertices, then $\gamma(G) + \Psi(G) = n$ where $\gamma(G)$ is the well-studied domination number of $G$. In this paper, we continue the study of the competition-enclaveless game introduced in 2001 by Philips and Slater [Graph Theory Notes N. Y. 41 (2001), 37--41] and defined as follows. Two players take turns in constructing a maximal enclaveless set $S$, where one player, Maximizer, tries to maximize $|S|$ and one player, Minimizer, tries to minimize~$|S|$. The competition-enclaveless game number $\Psi_g^+(G)$ of $G$ is the number of vertices played when Maximizer starts the game and both players play optimally. We study among other problems the conjecture that if $G$ is an isolate-free graph of order $n$, then $\Psi_g^+(G) \ge \frac{1}{2}n$. We prove this conjecture for regular graphs and for claw-free graphs.
\end{abstract}


\newpage
\section{Introduction}

A \emph{neighbor} of a vertex $v$ in $G$ is a vertex that is adjacent to $v$. A vertex \emph{dominates} itself and its neighbors. A \emph{dominating set} of a graph $G$ is a set $S$ of vertices of $G$ such that every vertex in $G$ is dominated by a vertex in $S$. The \emph{domination number} of $G$, denoted $\gamma(G)$, is the minimum cardinality of a dominating set in $G$, while the \emph{upper domination number} of $G$, denoted $\Gamma(G)$, is the maximum cardinality of a minimal dominating set in $G$. A minimal dominating set of cardinality~$\Gamma(G)$ we call a $\Gamma$-\emph{set of $G$}.

The \emph{open neighborhood} of a vertex $v$ in $G$ is the set of neighbors of $v$, denoted $N_G(v)$. Thus, $N_G(v) = \{u \in V \, | \, uv \in E(G)\}$. The \emph{closed neighborhood of $v$} is the set $N_G[v] = \{v\} \cup N_G(v)$. If the graph $G$ is clear from context, we simply write $N(v)$ and $N[v]$ rather than $N_G(v)$ and $N_G[v]$, respectively.

As defined by Alan Goldman and introduced in Slater~\cite{Sl77}, for a subset $S$ of vertices in a graph $G$, a vertex $v \in S$ is an \emph{enclave} of $S$ if it and all of its neighbors are also in $S$; that is, if $N[v] \subseteq S$. A set $S$ is \emph{enclaveless} if it does not contain any enclaves. We note that a set $S$ is a dominating set of a graph $G$ if the set $V(G) \setminus S$ is enclaveless. The \emph{enclaveless number} of $G$, denoted $\Psi(G)$, is the maximum cardinality of an enclaveless set in $G$, and the \emph{lower enclaveless number} of $G$, denoted by $\psi(G)$, is the minimum cardinality of a maximal enclaveless set. The domination and enclaveless numbers of a graph $G$ are related by the following equations.

\begin{obser}
\label{ob:relate}
If $G$ is a graph of order~$n$, then $\gamma(G) + \Psi(G) = n = \Gamma(G) + \psi(G)$.
\end{obser}

The domination game on a graph $G$ consists of two players, \emph{Dominator} and \emph{Staller}, who take turns choosing a vertex from $G$. Each vertex chosen must dominate at least one vertex not dominated by the vertices previously chosen. Upon completion of the game, the set of chosen (played) vertices is a dominating set in $G$. The goal of Dominator is to end the game with a minimum number of vertices chosen, while Staller has the opposite goal and wishes to end the game with as many vertices chosen as possible.

The Dominator-start domination game and the Staller-start domination game is the domination game when Dominator and Staller, respectively, choose the first vertex.  We refer to these simply as the D-game and S-game, respectively. The \emph{D-game domination number}, $\dstart(G)$, of $G$ is the minimum possible number of moves in a D-game when both players play optimally. The \emph{S-game domination number}, $\sstart(G)$, of $G$ is defined analogously for the S-game. The domination game was introduced by Bre{\v{s}}ar, Klav{\v{z}}ar, and Rall~\cite{BrKlRa10} and has been subsequently extensively studied in the literature (see, for example,~\cite{Bu2015a,Bu2015b,HeKi14,HeLo17,KiWeZa13,ko-2017,Sc17}).

Philips and Slater~\cite{PhSl01,PhSl02}
introduced what they called the \emph{competition}-\emph{enclaveless game}. The game is played by two players, Maximizer and Minimizer, on some graph $G$. They take turns in constructing a maximal enclaveless set $S$ of $G$. That is, in each turn a player plays a vertex $v$ that is not in the set $S$ of the vertices already chosen and such that $S \cup \{v\}$ does not contain an enclave, until there is no such vertex. We call such a vertex a \emph{playable vertex}. The goal of Maximizer is to make the final set $S$ as large as possible and for Minimizer to make the final set $S$ as small as possible.

The \emph{competition}-\emph{enclaveless game number}, or simply the  \emph{enclaveless game number}, $\Psi_g^+(G)$ of $G$ is the number of vertices chosen when Maximizer starts the game and both players play an optimal strategy according to the rules. The \emph{Minimizer-start competition}-\emph{enclaveless game number}, or simply the \emph{Minimizer-start enclaveless game number}, $\Psi_g^-(G)$, of $G$ is the number of vertices chosen when Minimizer starts the game and both players play an optimal strategy according to the rules. The competition-enclaveless game, which has been studied for example in~\cite{GoHe18,He18,PhSl01,PhSl02,SeSl07}, has not yet been explored in as much depth as the domination game. In this paper we continue the study of the competition-enclaveless game. Our main motivation for our study are the following conjectures that have yet to be settled, where an isolate-free graph is a graph that does not contain an isolated vertex.

\begin{conj}
\label{conj1}
If $G$ is an isolate-free graph of order $n$, then $\Psi_g^+(G) \ge \frac{1}{2}n$.
\end{conj}

Conjecture~\ref{conj1} was first posed as a question by Slater~\cite{Slater} to the 2nd author on 8th May 2015, and subsequently posed as a conjecture in~\cite{He18}. We refer to Conjecture~\ref{conj1} for general isolate-free graphs as the $\mathbf{\frac{1}{2}}$-\textbf{Enclaveless Game Conjecture}. We also pose the following conjecture for the Minimizer-start enclaveless game, where $\delta(G)$ denotes the minimum degree of the graph $G$.

\begin{conj}
\label{conj3}
If $G$ is a graph of order $n$ with $\delta(G) \ge 2$,
then $\Psi_g^-(G) \ge \frac{1}{2}n$.
\end{conj}

We proceed as follows. In Section~\ref{S:compare}, we discuss the domination game versus the enclaveless game, and show that these two games are very different and are not related. In Section~\ref{S:Fbounds}, we present fundamental bounds on the enclaveless game number and the Minimizer-start enclaveless game number. In Sections~\ref{S:regular} and~\ref{S:clawfree}, we show that the $\frac{1}{2}$-Enclaveless Game Conjecture holds for regular graphs and claw-free graphs, respectively. We use the standard notation $[k] = \{1,\ldots,k\}$.

\section{Game domination versus enclaveless game}
\label{S:compare}

Although the domination and enclaveless numbers of a graph $G$ are related by the equation $\gamma(G) + \Psi(G) = n$ (see Observation~\ref{ob:relate}), as remarked in~\cite{He18} the competition-enclaveless game is very different to the domination game. For example, if $k \ge 3$ and $G$ is a tree with exactly two non-leaf vertices both of which have $k$ leaf neighbors, that is, if $G$ is a double star $S(k,k)$, then $\Psi_g^+(G) = \Psi_g^-(G) = k+1$ and $\dstart(G) = 3$ and $\sstart(G) = 4$. If $n \ge 1$, then Ko\v{s}mrlj~\cite{ko-2017} showed that $\sstart(P_n) = \left\lceil \frac{n}{2} \right\rceil$ and that $\dstart(P_n) = \left\lceil\frac{n}{2}\right\rceil-1$ if $n \equiv 3 \,  (\modo \, 4)$ and $\dstart(P_n) = \left\lceil \frac{n}{2} \right\rceil$, otherwise. This is in contrast to the enclaveless game numbers of a path $P_n$ on $n \ge 2$ vertices determined by Phillips and Slater~\cite{PhSl02}.

\begin{theorem}{\rm (\cite{PhSl02})}
\label{enclave-path}
If $n \ge 2$, then $\Psi_g^+(P_n) = \lfloor  \frac{3n+1}{5}  \rfloor$ and $\Psi_g^-(P_n) = \lfloor  \frac{3n}{5} \rfloor$.
\end{theorem}

We remark that for the competition-enclaveless game the numbers $\Psi_g^+(G)$ and $\Psi_g^-(G)$ can vary greatly. For example, if $n \ge 1$ and $G$ is a star $K_{1,n}$, then $\Psi_g^+(G) = n$ while $\Psi_g^-(G) = 1$. However, for the domination game the Dominator-start game domination number and the Staller-start game domination number can differ by at most~$1$. The most significant difference between the domination game and the competition-enclaveless game is that the so-called Continuation Principle holds for the domination game but does not hold for the competition-enclaveless game.

Another significant difference between the domination game and the competition-enclaveless game is that upon completion of the domination game, the set of played vertices is a dominating set although not necessarily a minimal dominating set, while upon completion of the competition-enclaveless game, the set of played vertices is always a maximal enclaveless set. Thus, the enclaveless game numbers of a graph $G$ are always squeezed between the lower enclaveless number $\psi(G)$ of $G$ and the enclaveless number $\Psi(G)$ of $G$. We state this formally as follows.

\begin{obser}
\label{ob:bound}
If $G$ is a graph of order~$n$, then
\[
\psi(G) \le \Psi_g^-(G) \le \Psi(G) \hspace*{0.5cm} \mbox{and} \hspace*{0.5cm} \psi(G) \le \Psi_g^+(G) \le \Psi(G).
\]
\end{obser}

A graph $G$ is \emph{well}-\emph{dominated} if all the minimal dominating sets of $G$ have the same cardinality. Examples of well-dominated graphs include, for example, the complete graph $K_n$, $C_7$, $P_{10}$, the corona of any graph, and the graph formed from two vertex disjoint cycles of order~$5$ joined by a single edge. Finbow, Hartnell and Nowakowski~\cite{FiHaNo88} characterized the well-dominated graphs having no $3$-cycle nor $4$-cycle. As observed earlier, upon completion of the enclaveless game, the set of played vertices is always a maximal enclaveless set. Hence, any sequence of legal moves by Maximizer and Minimizer (regardless of strategy) in the enclaveless game will always lead to the game played on a graph $G$ of order~$n$ ending in $n - \gamma(G)$ moves. Thus as a consequence of Observation~\ref{ob:bound}, we have the following interesting connection between the enclaveless game and the class of well-dominated graphs.

\begin{obser}
\label{ob:well-dom}
If $G$ is a well-dominated graph of order~$n$, then $\Psi_g^-(G) = \Psi_g^+(G) = n - \gamma(G)$.
\end{obser}

It is well-known that if $G$ is an isolate-free graph of order $n$, then $\gamma(G) \le \frac{1}{2}n$, implying by Observation~\ref{ob:relate} that $\Psi(G) = n - \gamma(G) \ge \frac{1}{2}n$. Hence one might think that $\gamma_g(G) \le \Psi_g^+(G)$ for such a graph $G$ with no isolated vertex.  We now provide an infinite class of graphs to show that the ratio  $\gamma_g/\Psi_g^+$ of these two graphical invariants can be strictly larger than, and bounded away from,~$1$. The \emph{corona} $\coro(G)$ of a graph $G$, also denoted $G \circ K_1$ in the literature, is the graph obtained from $G$ by adding for each vertex $v$ of $G$ a new vertex $v'$ and the edge $vv'$ (and so, the vertex $v'$ has degree~$1$ in $\coro(G)$). The edge $vv'$ is called a \emph{pendant edge}.

\begin{theorem}
\label{t:ratio}
If $n \ge 2$ is an integer and $\cG_n$ denotes the class of all isolate-free graphs $G$ of order~$n$, then
\[
\sup_{n} \, \frac{ \gamma_g(G) }{\Psi_g^+(G)} \ge \frac{11}{10}
\]
where the supremum is taken over all graphs $G \in \cG_n$.
\end{theorem}
\proof Let $n=10q$ for some positive integer $q$ and $G_n$ be the corona of the path $P_n$.  That is, the vertex set of $G_n$ is $X_n \cup Y_n$ where $X_n = \{x_i : i \in [n]\}$ and $Y_n=\{y_i : i \in [n]\}$.  The edge set of $G_n$ is $\{x_ix_{i+1} : i \in [n-1] \} \cup \{x_iy_i : i \in [n]\}$.  For each $k$ such that $0 \le k \le q-1$ we let $B_k$ be the subgraph of $G_n$ induced by $\cup_{i=1}^{10} \{x_{10k+i}, y_{10k+i}\}$.  The D-game is played on $G_n$.  At any point in this game we say that $B_i$ is open if no vertex in $B_i$ has been played by either player; otherwise we say $B_i$ is not open.  By the Continuation Principle we may assume that any vertex played by Dominator belongs to $X_n$. We denote by $d_1, d_2, \ldots $ and  $s_1, s_2, \ldots $ the sequence of moves played by Dominator and Staller in the domination game.  We now provide a strategy for Staller to show that $\gamma_g(G_n) \ge 11q$.

\begin{enumerate}
\item If Dominator plays $d_1=x_i$ where $10k+1 \le i \le 10k+5$ for some $k$ such that $0 \le k \le q-1$, then Staller plays $s_1=y_{10k+8}$.
\item If Dominator plays $d_1=x_i$ where $10k+6 \le i \le 10k+10$ for some $k$ such that $0 \le k \le q-1$, then Staller plays $s_1=y_{10k+3}$.
\end{enumerate}

 If Dominator plays a vertex $d_j$ in an open $B_i$ for some $j \ge 1$ and $i$ with $0 \le i \le q-1$, then Staller plays $s_j$ also in $B_i$ as described in (a) and (b) above.  On the other hand, suppose that Dominator plays a vertex $d_j$ in a $B_i$ that is not open. If this move of Dominator is his second move played in $B_i$, then in this case Staller plays the support vertex that is adjacent to the leaf she played earlier in the game when $B_i$ changed from being open to being not open. This support vertex was a legal move for Staller because of the structure of the graph $G_n$. When the game ends at least one vertex for each pair $x_k,y_k$ must have been played by one of the players.  Furthermore, the above strategy for Staller shows that she can ensure that at least eleven vertices are played from each of $B_0, \ldots, B_{q-1}$.  Therefore, $\gamma_g(G_n) \ge 11q$.

Now, observe that every minimal dominating set of $G_n$ has cardinality $n$, which implies by Observation~\ref{ob:relate} that every maximal enclaveless set of $G_n$ also has cardinality $n$; that is, $\psi(G) = \Psi(G) = n$ where we recall that $\psi(G)$ denotes the cardinality of the smallest maximal enclaveless set in $G$ and $\Psi(G)$ is the
cardinality of a largest enclaveless set in $G$.  Hence by Observation~\ref{ob:bound}, $\Psi_g^+(G_n) = n$. Consequently, we have shown
\[
\sup_n \frac{\gamma_g(G_n)}{\Psi_g^+(G_n)} \ge \frac{11}{10}\,.
\]
The desired result follows noting that $G_n \in \cG_{2n}$.~\QED

\section{Fundamental bounds}
\label{S:Fbounds}

In this section, we establish some fundamental bounds on the (Maximizer-start) enclaveless game number and the Minimizer-start enclaveless game number. We establish next an upper bound on the enclaveless number of a graph in terms of the maximum degree and order of the graph.

\begin{prop}
\label{p:bound1}
If $G$ is an isolate-free  graph of order~$n$ with maximum degree $\Delta(G) = \Delta$, then
\[
\left( \frac{1}{\Delta+1} \right) n \le \psi(G) \le \Psi(G) \le \left( \frac{\Delta}{\Delta+1} \right) n.
\]
\end{prop}
\proof If $G$ is any isolate-free graph of order $n$  and maximum degree $\Delta$, then $\gamma(G) \ge \frac{n}{\Delta + 1}$, with equality precisely when $G$ has a minimum dominating set consisting of vertices of degree~$\Delta$ that is a $2$-packing, where a $2$-packing is a set $S$ of vertices that are pairwise at distance at least~$3$ apart. Hence, by Observation~\ref{ob:relate},
\[
\Psi(G) = n - \gamma(G) \le n - \frac{n}{\Delta+1}= \left( \frac{\Delta}{\Delta+1} \right) n.
\]

On the other hand, let $D$ be a minimal dominating set of maximum cardinality, and so $|D| = \Gamma(G)$. Let $\barD = V(G) \setminus D$, and so $|\barD| = n - |D|$. Let $\ell$ be the number of edges between $D$ and $\barD$. Since $D$ is a minimal dominating set, every vertex in $D$ has at least one neighbor in $\barD$, and so $\ell \ge |D|$. Since $G$ has maximum degree~$\Delta$, every vertex in $\barD$ has at most $\Delta$ neighbors in $D$, and so $\ell \le
\Delta \cdot |\barD| = \Delta (n - |D|)$. Hence, $|D| \le \Delta (n - |D|)$, implying that $\Gamma(G) = |D| \le \Delta n /(\Delta+1)$. Thus by Observation~\ref{ob:relate},
\[
\psi(G) = n - \Gamma(G) \ge n - \left( \frac{\Delta}{\Delta+1} \right) n = \left( \frac{1}{\Delta+1} \right) n.
\]
This completes the proof of Proposition~\ref{p:bound1}.~\QED

\medskip
By Observation~\ref{ob:bound}, the set of played vertices in either the Maximizer-start enclaveless game or the Minimizer-start enclaveless game is an
enclaveless set of $G$. Thus as an immediate consequence of Proposition~\ref{p:bound1}, we have the following result.

\begin{prop}
\label{prop:trivialupperbound}
If $G$ is an isolate-free graph of order~$n$ with maximum degree $\Delta(G) = \Delta$, then
\[
\left( \frac{1}{\Delta+1} \right) n \le \Psi_g^-(G) \le \left( \frac{\Delta}{\Delta+1} \right) n \hspace*{0.5cm} \mbox{and} \hspace*{0.5cm}
\left( \frac{1}{\Delta+1} \right) n \le \Psi_g^+(G) \le \left( \frac{\Delta}{\Delta+1} \right) n.
\]
\end{prop}

\medskip
The lower bound in Proposition~\ref{prop:trivialupperbound} on $\Psi_g^-(G)$ is achieved, for example, by taking $G = K_{1,\Delta}$ for any given $\Delta \ge 1$ in which case $\Psi_g^-(G) = 1 = ( \frac{1}{\Delta+1} ) n$ where $n = n(G) = \Delta + 1$. We show next that the upper bounds in Proposition~\ref{prop:trivialupperbound} are realized for infinitely many connected graphs.

\begin{prop}
\label{prop:achieve1}
There exist infinitely many positive integers $n$ along with a connected graph $G$ of order~$n$ satisfying
\[
\Psi_g^-(G) = \Psi_g^+(G) = \left( \frac{\Delta(G)}{\Delta(G)+1} \right) n.
\]
\end{prop}
\proof Let $r$ be an integer such that $r \ge 4$ and let $m$ be any positive integer.  For each $i \in [m]$, let $H_i$ be a graph obtained from a complete graph of order $r+1$ by removing the edge $x_iy_i$ for two distinguished vertices $x_i$ and $y_i$. The graph $F_m$ is obtained from the disjoint union of $H_1,\ldots,H_m$ by adding the edges $y_ix_{i+1}$ for each $i \in [m]$ where the subscripts are computed modulo $m$.  The vertices $x_i$ and $y_i$ are called connectors in $F_m$, and each of the $r-1$ vertices in the set $V(H_i) \setminus \{x_i,y_i\}$ is called a hidden vertex of $H_i$.  Note that $F_m$ is $r$-regular and has order $n=m(r+1)$.

We first show  that $\Psi_g^-(F_m) = ( \frac{r}{r+1} ) n$. Suppose the Minimizer-start enclaveless game is played on $F_m$.  We provide a strategy for Maximizer that forces exactly $rm$ vertices to be played. Maximizer's strategy is to make sure that all the connector vertices in the graph are played.  If he can accomplish this, then exactly $rm$ vertices will be played when the game ends because of the structure of $F_m$.  Suppose that at some point in the game Minimizer plays a vertex from some $H_j$.  If one of the connector vertices, say $x_j$, is playable, then Maximizer responds by playing $x_j$.  If both connector vertices have already been played and some hidden vertex, say $w$, in $H_j$ is playable, then Maximizer plays $w$.  If no vertex of $H_j$ is playable, then Maximizer plays a connector vertex from $H_i$ for some $i \neq j$ if one is playable and otherwise plays any playable vertex.  Since $H_k$ contains at least $3$ hidden vertices for each $k \in [m]$, it follows that Maximizer can guarantee that all the connector vertices are played by following this strategy. This implies that for each $i \in [m]$, exactly one hidden vertex of $H_i$ is not played during the course of the game.  That is, the set of played vertices has cardinality
\[
rm = \left( \frac{r}{r+1} \right) m(r+1) = \left( \frac{\Delta(F_m)}{\Delta(F_m)+1} \right) n\,,
\]
where we recall that $\Delta(F_m) = r$. Thus,
\[
\Psi_g^-(F_m) \ge \left( \frac{\Delta(F_m)}{\Delta(F_m)+1} \right) n.
\]
By Proposition~\ref{prop:trivialupperbound}, \[
\Psi_g^-((F_m)) \le \left( \frac{\Delta(F_m)}{\Delta(F_m)+1} \right) n.
\]
Consequently, $\Psi_g^-(F_m) = ( \frac{\Delta(F_m)}{\Delta(F_m)+1} ) n$.

If the Maximizer-start enclaveless game is played on $F_m$, then the same strategy as above for Maximizer forces $rm$ vertices to be played (even with the relaxed condition that $r$ be an integer larger than $2$). Thus as before, $\Psi_g^+(F_m) = ( \frac{\Delta(F_m)}{\Delta(F_m)+1} ) n$.~\QED

\section{Regular graphs}
\label{S:regular}

In this section, we show that $\frac{1}{2}$-Enclaveless Game Conjecture (see Conjecture~\ref{conj1}) holds for the class of regular graphs, as does Conjecture~\ref{conj3} for the Minimizer-start enclaveless game. For a set $S \subset V(G)$ of vertices in a graph $G$ and a vertex $v \in S$, we define the \emph{$S$-external private neighborhood} of a vertex $v$, abbreviated $\epn_G(v,S)$, as the set of all vertices outside $S$ that are adjacent to $v$ but to no other vertex of $S$; that is,
\[
\epn_G(v,S) = \{w \in V(G) \setminus S \mid N_G(w) \cap S = \{v\}\}.
\]
We define an \emph{$S$-external private neighbor} of $v$ to be a vertex in $\epn_G(v,S)$.

\begin{theorem}\label{t:regular}
If $G$ is a $k$-regular graph of order $n$, then $\Psi_g^+(G) \ge \frac{1}{2}n$ and $\Psi_g^-(G) \ge \frac{1}{2}n$.
\end{theorem}
\proof
Suppose the Maximizer-start enclaveless game is played on $G$.  Let $S$ denote the set of all vertices played when
the game ends.  By definition of the game, the set $S$ is a maximal enclaveless set in $G$. By Observations~\ref{ob:relate} and~\ref{ob:bound}, we have $|S| = \Psi_g^+(G) \ge \psi(G) = n-\Gamma(G)$. It therefore suffices to establish the proposition by proving that $\Gamma(G) \le \frac{1}{2}n$.

This inequality is proved in~\cite{SoHe13}, but we prove it here for the sake of completeness. Let $D$ be an arbitrary minimal dominating set of $G$.  Denote by $D_1$ the set of vertices in $D$ that have a $D$-external private neighbor.  That is, $D_1 = \{x \in D : \epn_G(x,D) \ne \emptyset\}$.  In addition, let $D_2 = D \setminus D_1$.  Since $D$ is a minimal dominating set, the set $D_2$ consists of those vertices in $D$ that are isolated in the subgraph $G[D]$ of $G$ induced by $D$.  Let
\[
C_1=\bigcup_{x\in D_1} \epn_G(x,D) \hspace*{0.5cm} \mbox{and} \hspace*{0.5cm} C_2 = V(G) \setminus (D \cup C_1).
\]

We note that by definition, there are no edges in $G$ joining a vertex in $D_2$ and a vertex in $C_1$. That is, each vertex in $D_2$ has $k$ neighbors in $C_2$.  Since every vertex has degree $k$, each vertex of $C_2$ has at most $k$
neighbors in $D_2$.  Denote by $\ell$ the number of edges of the form $uv$ where $u\in D_2$ and $v \in C_2$.
It now follows that $k|D_2|=\ell \le k|C_2|$.  That is, $|D_2| \le |C_2|$.  Now
\[
|D|=|D_1|+|D_2| \le |C_1|+|C_2|=n-|D|\,,
\]
which shows that $\Gamma(G)\le |D| \le \frac{1}{2}n$. Similarly, $\Psi_g^-(G) \ge \psi(G) = n-\Gamma(G) \ge \frac{1}{2}n$.~\QED

\medskip
We remark that the lower bound in Theorem~\ref{t:regular} is achieved for $k=1$ and $k=2$ as shown by $K_2$ and $C_4$, respectively. However, it remains an open problem to characterize the graphs achieving equality in Theorem~\ref{t:regular} for each value of $k \ge 1$.

A similar proof to that of Theorem~\ref{t:regular} will establish the same lower bounds by restricting the minimum degree and forbidding induced stars of a certain size.

\begin{prop} \label{prop:nolargestars}
Let $k$ be a positive integer.  If $G$ is a graph of order $n$ with minimum degree at least $k$ and with no induced $K_{1,k+1}$, then both $\Psi_g^+(G)$ and $\Psi_g^-(G)$ are at least $\frac{1}{2}n$.
\end{prop}
\proof
Let $D$ be a minimal dominating set of $G$.   The sets $D_1,D_2,C_1$ and $C_2$  as well as $\ell$ are defined as in the proof of Theorem~\ref{t:regular}.  In this case we get $k|D_2| \le \ell$ and $\ell \le k|C_2|$.  The first of these inequalities follows since $\delta(G) \ge k$ and the second inequality follows from the fact that $D_2$ is independent and the assumption that $G$ is $K_{1,k+1}$-free.  Once again we conclude  that $|D_2| \le |C_2|$, and the result follows.~\QED

\section{Claw-free graphs}
\label{S:clawfree}

A graph is \emph{claw}-\emph{free} if it does not contain the star $K_{1,3}$ as an induced subgraph. In this section, we show that $\frac{1}{2}$-Enclaveless Game Conjecture (see Conjecture~\ref{conj1}) holds for the class of claw-free graphs with no isolated vertex, as does Conjecture~\ref{conj3} for the Minimizer-start enclaveless game.
For this purpose, we recall the definition of an irredundant set. For a set $S$ of vertices in a graph $G$ and a vertex $v \in S$, the \emph{$S$-private neighborhood} of $v$ is the set
\[
\pn_G[v,S] = \{w \in V \mid N_G[w] \cap S = \{v\}\}.
\]

If the graph $G$ is clear from context, we simply write $\pn[v,S]$ rather than $\pn_G[v,S]$. We note that if the vertex $v$ is isolated in $G[S]$, then $v \in \pn[v,S]$. A vertex in the set $\pn[v,S]$ is called an $S$-\emph{private neighbor} of $v$. The set $S$ is an \emph{irredundant set} if every vertex of $S$ has an $S$-private neighbor. The \emph{upper irredundance number} $\IR(G)$ is the maximum cardinality of an irredundant set in $G$.

The \emph{independence number} $\alpha(G)$ of $G$ is the maximal cardinality of an independent set of vertices in $G$. An independent set of vertices of $G$ of cardinality $\alpha(G)$ we call an $\alpha$-\emph{set of $G$}. Every maximum independent set in a graph is minimal dominating, and every minimal dominating set is irredundant. Hence we have the following inequality chain.

\begin{obser}{\rm (\cite{CoHeMi78})}
\label{ob:dom_chain}
For every graph $G$, we have $\alpha(G) \le \Gamma(G) \le \IR(G)$.
\end{obser}

The inequality chain in Observation~\ref{ob:dom_chain} is part of the canonical domination chain which was first observed by Cockayne, Hedetniemi, and Miller~\cite{CoHeMi78} in 1978. We shall need the following upper bounds on the independence number of a claw-free graph.

\begin{theorem}
\label{t:indep}
If $G$ is a connected claw-free graph of order~$n$ and minimum degree~$\delta \ge 1$, then the following holds.
\\ [-26pt]
\begin{enumerate}
\item  {\rm (\cite{Ga99,RySc95})} If $\delta = 1$, then $\alpha(G) \le \frac{1}{2}(n+1)$. \2
\item {\rm (\cite{Faudree92,LiVi90})} If $\delta \ge 2$, then $\alpha(G) \le \frac{2n}{\delta + 2}$.
\end{enumerate}
\end{theorem}

In 2004, Favaron~\cite{Fa03} established the following upper bound on the irredundance number of a claw-free graph.

\begin{theorem}{\rm (\cite{Fa03})}
\label{t:bound_IR1}
If $G$ is a connected, claw-free graph of order~$n$, then $\IR(G) \le \frac{1}{2}(n+1)$. Moreover, if $\IR(G) = \frac{1}{2}(n+1)$, then $\alpha(G) = \Gamma(G) = \IR(G)$.
\end{theorem}

If $G$ is a connected, claw-free graph of order~$n$ and minimum degree~$\delta \ge 2$, then by Theorem~\ref{t:indep}(b) we have $\alpha(G) \le \frac{1}{2}n$. In this case when $\delta \ge 2$, if $\IR(G) = \frac{1}{2}(n+1)$ holds, then by Theorem~\ref{t:bound_IR1} we have $\alpha(G) = \frac{1}{2}(n+1)$, a contradiction. Hence when $\delta \ge 2$, we must have $\IR(G) \le \frac{1}{2}n$. We state this formally as follows.

\begin{cor}{\rm (\cite{Fa03})}
\label{c:bound_IR1}
If $G$ is a connected, claw-free graph of order~$n$ and minimum degree at least~$2$, then $\IR(G) \le \frac{1}{2}n$.
\end{cor}

We are now in a position to prove the following result.

\begin{theorem}\label{t:clawfree1}
If $G$ is a connected claw-free graph of order $n$ and $\delta(G) \ge 2$, then
\[
\Psi_g^+(G) \ge \frac{1}{2}n \hspace{0.5cm} \mbox{and} \hspace{0.5cm} \Psi_g^-(G) \ge \frac{1}{2}n.
\]
\end{theorem}
\proof
Suppose the Minimizer-start enclaveless game is played on $G$. Let $S$ denote the set of all vertices played when
the game ends.  By definition of the game, the set $S$ is a maximal enclaveless set in $G$. By Observations~\ref{ob:relate},~\ref{ob:bound} and~\ref{ob:dom_chain} and Corollary~\ref{c:bound_IR1}, we have
\[
|S| = \Psi_g^-(G) \ge \psi(G) = n - \Gamma(G) \ge n - \IR(G) \ge n - \frac{1}{2}n = \frac{1}{2}n,
\]
as desired. Similarly, $\Psi_g^+(G) \ge \psi(G) \ge n - \IR(G) \ge \frac{1}{2}n$.~\QED

\medskip
By Theorem~\ref{t:clawfree1}, we note that Conjecture~\ref{conj3} holds for connected claw-free graphs. In order to prove that Conjecture~\ref{conj1} holds for connected claw-free graphs, we need the characterization due to Favaron~\cite{Fa03} of the graphs achieving equality in the bound of Theorem~\ref{t:bound_IR1}. For this purpose, we recall that a vertex $v$ of a graph $G$ is a \emph{simplicial vertex} it is neighborhood $N_G(v)$ induces a complete subgraph of $G$. A \emph{clique} of a graph $G$ is a maximal complete subgraph of $G$. The \emph{clique graph} of $G$ has the set of cliques of $G$ as its vertex set, and two vertices in the clique graph are adjacent if and only if they intersect as cliques of $G$. A \emph{non}-\emph{trivial tree} is a tree of order at least~$2$.

Favaron~\cite{Fa03} defined the family $\cF$ of claw-free graphs $G$ as follows. Let $T_1, \ldots, T_q$ be $q \ge 1$ non-trivial trees. Let $L_i$ be the line graph of the corona $\coro(T_i)$ of the tree $T_i$ for $i \in [q]$. If $q = 1$, let $G = L_1$. If $q \ge 2$, let $G$ be the graph constructed from the line graphs $L_1, L_2, \ldots, L_q$ by choosing $q-1$ pairs $\{x_{ij},x_{ji}\}$ such that the following holds. \\ [-24pt]
\begin{enumerate}
\item[$\bullet$] $x_{ij}$ and $x_{ji}$ are simplicial vertices of $L_i$ and $L_j$, respectively, where $i \ne j$.
\item[$\bullet$] The $2(q-1)$ vertices from the $q-1$ pairs $\{x_{ij},x_{ji}\}$ are all distinct vertices.
\item[$\bullet$] Contracting each pair of vertices $x_{ij}$ and $x_{ji}$ into one common vertex $c_{ij}$ results in a graph whose clique graph is a tree.
\end{enumerate}

\begin{figure}[htb]
\begin{center}
\begin{tikzpicture}[scale=.7,style=thick,x=1cm,y=1cm]
\def\vr{2.75pt}
\def\vB{14pt}
\path (1,0) coordinate (v1);
\path (2,1) coordinate (v2);
\path (3,0) coordinate (v3);
\path (4,1) coordinate (v4);
\path (5,0) coordinate (v5);
\path (6,1) coordinate (v6);
\path (6,-0.15) coordinate (v6p);
\path (7,0) coordinate (v7);
\path (8,1) coordinate (v8);
\path (9,0) coordinate (v9);
\path (10,1) coordinate (v10);
\path (11,0) coordinate (v11);
\path (12,1) coordinate (v12);
\path (13,0) coordinate (v13);
\draw (v1) -- (v2);
\draw (v2) -- (v3);
\draw (v3) -- (v4);
\draw (v4) -- (v5);
\draw (v5) -- (v6);
\draw (v5) -- (v8);
\draw (v6) -- (v7);
\draw (v8) -- (v9);
\draw (v9) -- (v10);
\draw (v10) -- (v11);
\draw (v11) -- (v12);
\draw (v12) -- (v13);
\draw (v4) -- (v6);
\draw (v6) -- (v8);
\draw (v10) -- (v12);
\draw (v1) [fill=white] circle (\vr);
\draw (v2) [fill=black] circle (\vr);
\draw (v3) [fill=white] circle (\vr);
\draw (v4) [fill=black] circle (\vr);
\draw (v5) [fill=white] circle (\vr);
\draw (v6) [fill=black] circle (\vr);
\draw (v7) [fill=white] circle (\vr);
\draw (v8) [fill=black] circle (\vr);
\draw (v9) [fill=white] circle (\vr);
\draw (v10) [fill=black] circle (\vr);
\draw (v11) [fill=white] circle (\vr);
\draw (v12) [fill=black] circle (\vr);
\draw (v13) [fill=white] circle (\vr);
\draw (v4) to[out=60,in=120, distance=1cm] (v8);
\draw[anchor = north] (v3) node {{\small $c_{12}$}};
\draw[anchor = north] (v9) node {{\small $c_{23}$}};
\draw[anchor = north] (v6p) node {{\small $G$}};
\path (-2,2.5) coordinate (u0);
\path (-1,3) coordinate (u1);
\path (0,4) coordinate (u2);
\path (0,2.85) coordinate (u2p);
\path (1,3) coordinate (u3);
\path (3,3) coordinate (u4);
\path (4,4) coordinate (u5);
\path (5,3) coordinate (u6);
\path (6,4) coordinate (u7);
\path (6,2.85) coordinate (u7p);
\path (7,3) coordinate (u8);
\path (8,4) coordinate (u9);
\path (9,3) coordinate (u10);
\path (11,3) coordinate (u11);
\path (12,4) coordinate (u12);
\path (13,3) coordinate (u13);
\path (13,2.85) coordinate (u13p);
\path (14,4) coordinate (u14);
\path (15,3) coordinate (u15);
\draw (u1) -- (u2);
\draw (u2) -- (u3);
\draw (u4) -- (u5);
\draw (u5) -- (u6);
\draw (u6) -- (u7);
\draw (u6) -- (u9);
\draw (u7) -- (u8);
\draw (u9) -- (u10);
\draw (u11) -- (u12);
\draw (u12) -- (u13);
\draw (u13) -- (u14);
\draw (u14) -- (u15);
\draw (u5) -- (u7);
\draw (u7) -- (u9);
\draw (u12) -- (u14);
\draw (u1) [fill=white] circle (\vr);
\draw (u2) [fill=black] circle (\vr);
\draw (u3) [fill=white] circle (\vr);
\draw (u4) [fill=white] circle (\vr);
\draw (u5) [fill=black] circle (\vr);
\draw (u6) [fill=white] circle (\vr);
\draw (u7) [fill=black] circle (\vr);
\draw (u8) [fill=white] circle (\vr);
\draw (u9) [fill=black] circle (\vr);
\draw (u10) [fill=white] circle (\vr);
\draw (u11) [fill=white] circle (\vr);
\draw (u12) [fill=black] circle (\vr);
\draw (u13) [fill=white] circle (\vr);
\draw (u14) [fill=black] circle (\vr);
\draw (u15) [fill=white] circle (\vr);
\draw (u5) to[out=60,in=120, distance=1cm] (u9);
\draw[anchor = north] (u3) node {{\small $x_{12}$}};
\draw[anchor = north] (u4) node {{\small $x_{21}$}};
\draw[anchor = north] (u10) node {{\small $x_{23}$}};
\draw[anchor = north] (u11) node {{\small $x_{32}$}};
\draw[anchor = north] (u2p) node {{\small $L_1$}};
\draw[anchor = north] (u7p) node {{\small $L_2$}};
\draw[anchor = north] (u13p) node {{\small $L_3$}};
\draw[anchor = north] (u0) node {{\small $\Downarrow$}};
\path (-2,5.5) coordinate (w0);
\path (-1,6) coordinate (w1);
\path (-1,7) coordinate (w2);
\path (0,5.85) coordinate (w2p);
\path (1,7) coordinate (w3);
\path (1,6) coordinate (w4);
\path (3,6) coordinate (w5);
\path (3,7) coordinate (w6);
\path (5,6) coordinate (w7);
\path (5,7) coordinate (w8);
\path (6,5.85) coordinate (w8p);
\path (7,6) coordinate (w9);
\path (7,7) coordinate (w10);
\path (9,6) coordinate (w11);
\path (9,7) coordinate (w12);
\path (11,6) coordinate (w13);
\path (11,7) coordinate (w14);
\path (13,6) coordinate (w15);
\path (13,7) coordinate (w16);
\path (13,5.75) coordinate (w15p);
\path (15,6) coordinate (w17);
\path (15,7) coordinate (w18);
\draw (w1) -- (w2);
\draw (w2) -- (w3);
\draw (w3) -- (w4);
\draw (w5) -- (w6);
\draw (w7) -- (w8);
\draw (w9) -- (w10);
\draw (w11) -- (w12);
\draw (w6) -- (w8);
\draw (w8) -- (w10);
\draw (w13) -- (w14);
\draw (w15) -- (w16);
\draw (w17) -- (w18);
\draw (w14) -- (w16);
\draw (w16) -- (w18);
%
\draw (w1) [fill=white] circle (\vr);
\draw (w2) [fill=black] circle (\vr);
\draw (w3) [fill=black] circle (\vr);
\draw (w4) [fill=white] circle (\vr);
\draw (w5) [fill=white] circle (\vr);
\draw (w6) [fill=black] circle (\vr);
\draw (w7) [fill=white] circle (\vr);
\draw (w8) [fill=black] circle (\vr);
\draw (w9) [fill=white] circle (\vr);
\draw (w10) [fill=black] circle (\vr);
\draw (w11) [fill=white] circle (\vr);
\draw (w12) [fill=black] circle (\vr);
\draw (w13) [fill=white] circle (\vr);
\draw (w14) [fill=black] circle (\vr);
\draw (w15) [fill=white] circle (\vr);
\draw (w16) [fill=black] circle (\vr);
\draw (w17) [fill=white] circle (\vr);
\draw (w18) [fill=black] circle (\vr);
\draw (w8) to[out=60,in=120, distance=1cm] (w12);
\draw[anchor = north] (w2p) node {{\small $\coro(T_1)$}};
\draw[anchor = north] (w8p) node {{\small $\coro(T_2)$}};
\draw[anchor = north] (w15p) node {{\small $\coro(T_3)$}};
\draw[anchor = north] (w0) node {{\small $\Downarrow$}};
\path (-2,8.5) coordinate (x0);
\path (-1,9) coordinate (x1);
\path (1,9) coordinate (x2);
\path (0,8.85) coordinate (x2p);
\path (3,9) coordinate (x3);
\path (5,9) coordinate (x4);
\path (6,8.85) coordinate (x4p);
\path (7,9) coordinate (x5);
\path (9,9) coordinate (x6);
\path (11,9) coordinate (x7);
\path (13,9) coordinate (x8);
\path (13,8.75) coordinate (x8p);
\path (15,9) coordinate (x9);
\draw (x1) -- (x2);
\draw (x3) -- (x4);
\draw (x4) -- (x5);
\draw (x7) -- (x8);
\draw (x8) -- (x9);
\draw (x1) [fill=black] circle (\vr);
\draw (x2) [fill=black] circle (\vr);
\draw (x3) [fill=black] circle (\vr);
\draw (x4) [fill=black] circle (\vr);
\draw (x5) [fill=black] circle (\vr);
\draw (x6) [fill=black] circle (\vr);
\draw (x7) [fill=black] circle (\vr);
\draw (x8) [fill=black] circle (\vr);
\draw (x9) [fill=black] circle (\vr);
\draw (x4) to[out=60,in=120, distance=1cm] (x6);
\draw[anchor = north] (x2p) node {{\small $T_1$}};
\draw[anchor = north] (x4p) node {{\small $T_2$}};
\draw[anchor = north] (x8p) node {{\small $T_3$}};
\draw[anchor = north] (x0) node {{\small $\Downarrow$}};
\end{tikzpicture}
\end{center}
\vskip -0.6cm
\caption{An illustration of the construction of a graph $G$ in the family $\cF$} \label{f:familyF}
\end{figure}
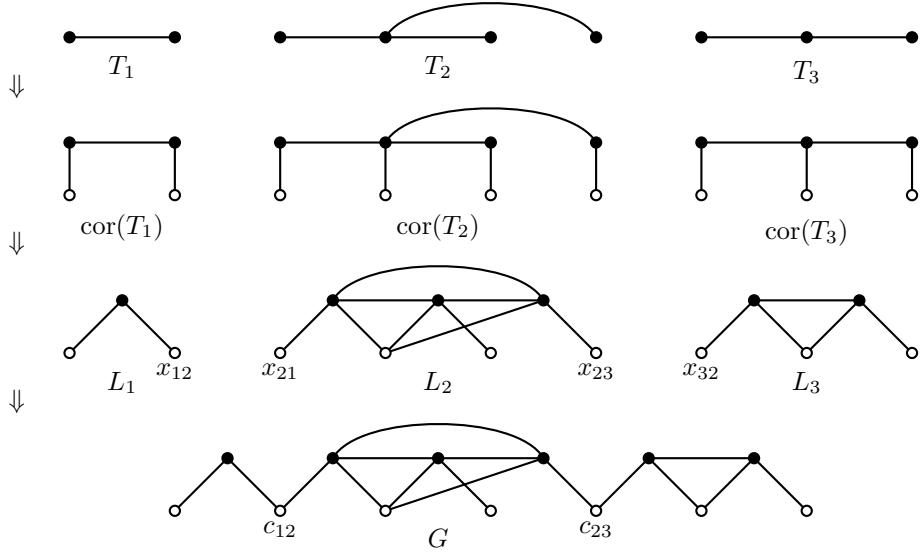

To illustrate the above construction of a graph $G$ in the family $\cF$ consider, for example, such a construction when $q = 3$ and the trees $T_1, T_2, T_3$ are given in Figure~\ref{f:familyF}.

We note that if $G$ is an arbitrary graph of order~$n$ in the family~$\cF$, then $n \ge 3$ is odd and the vertex set of $G$ can be partitioned into two sets $A$ and $B$ such that the following holds. \\ [-24pt]
\begin{enumerate}
\item[$\bullet$] $|A| = \frac{1}{2}(n-1)$ and $|B| = \frac{1}{2}(n+1)$.
\item[$\bullet$] The set $B$ is an independent set.
\item[$\bullet$] Each vertex in $A$ has exactly two neighbors in $B$.
\end{enumerate}

We refer to the partition $(A,B)$ as the partition associated with $G$. For the graph $G \in \cF$ illustrated in Figure~\ref{f:familyF}, the set $A$ consists of the darkened vertices and the set $B$ consists of the white vertices.

We are now in a position to state the characterization due to Favaron~\cite{Fa03} of the graphs achieving equality in the bound of Theorem~\ref{t:bound_IR1}.

\begin{theorem}{\rm (\cite{Fa03})}
\label{t:bound_IR2}
If $G$ is a connected, claw-free graph of order~$n \ge 3$, then $\IR(G) \le \frac{1}{2}(n+1)$, with equality if and only if $G \in \cF$.
\end{theorem}

We prove next the following property of graphs in the family~$\cF$.

\begin{lemma}
\label{l:lemma1}
If $G \in \cF$ and $(A,B)$ is the partition associated with $G$, then the set $B$ is the unique $\IR$-set of $G$.
\end{lemma}
\proof We proceed by induction on the order $n \ge 3$ of $G \in \cF$. If $n = 3$, then $G = P_3$. In this case, the set $B$ consists of the two leaves of $G$, and the desired result is immediate. This establishes the base case. Suppose that $n \ge 5$ and that the result holds for all graphs $G' \in \cF$ of order~$n'$, where $3 \le n' < n$. Let $Q$ be an $\IR$-set of $G$.

By construction of the graph $G$, the set $B$ contains at least two vertices of degree~$1$ in $G$. Let $v$ be an arbitrary vertex in $B$ of degree~$1$ in $G$, and let $u$ be its neighbor. We note that $u \in A$. Let $G' = G - \{u,v\}$ and let $G'$ have order~$n'$, and so $n' = n - 2$. Let $A' = A \setminus \{u\}$ and $B' = B \setminus \{v\}$.  By construction of the graph $G$ and our choice of the vertex~$v$, we note that $G' \in \cF$ and that $(A',B')$ is the partition associated with $G'$. Applying the inductive hypothesis to $G'$, the set $B'$ is the unique $\IR$-set of $G'$. Let $w$ be the second neighbor of $u$ in $G$ that belongs to the set $B$, and so $N_G(u) \cap B = \{v,w\}$. By the structure of the graph $G \in \cF$, we note that $N_G[w] \subset N_G[u]$ and that the subgraph of $G$ induced by $N_G[w]$ is a clique.

Suppose, to the contrary, that $Q \ne B$. Let $Q'$ be the restriction of $Q$ to the graph $G'$, and so $Q' = Q \cap V(G')$. Suppose that $u \in Q$. Since $Q$ is an irredundant set, this implies that $v \notin Q$. If $w \in Q$, then $\pn[w,Q] = \emptyset$, contradicting the fact that $Q$ is an irredundant set.  Hence, $w \notin Q$, and so $Q' \ne B'$. By the inductive hypothesis, the set $Q'$ is therefore not an $\IR$-set of $G'$, and so $|Q'| < \IR(G')$. Thus, $\IR(G) = |Q| = |Q'| + 1 \le (\IR(G') - 1) + 1 = \frac{1}{2}(n'+1) = \frac{1}{2}(n-1) < \IR(G)$, a contradiction. Hence, $u \notin Q$. In this case, $\IR(G) = |Q| \le |Q'| + 1 \le \IR(G') + 1 = \frac{1}{2}(n'+1) + 1 = \frac{1}{2}(n+1) = \IR(G)$. Hence, we must have equality throughout this inequality chain. This implies that $v \in Q$ and $|Q'| = \IR(G')$. By the inductive hypothesis, we therefore have $Q' = B'$. Hence, $Q = Q' \cup \{v\} = B' \cup \{v\} = B$. Thus, the set $B$ is the unique $\IR$-set of $G$.~\QED

\begin{cor}
\label{c:cor1}
If $G \in \cF$ and $(A,B)$ is the partition associated with $G$, then the set $B$ is the unique $\alpha$-set of $G$ and the unique $\Gamma$-set of $G$.
\end{cor}
\proof By Theorem~\ref{t:bound_IR1}, $\alpha(G) = \Gamma(G) = \IR(G) = \frac{1}{2}(n+1)$. By Lemma~\ref{l:lemma1}, the set $B$ is the unique $\IR$-set of $G$. Since every $\alpha$-set of $G$ is an $\IR$-set of $G$ and $\alpha(G) = \IR(G)$, this implies that $B$ is the unique $\alpha$-set of $G$. Since every $\Gamma$-set of $G$ is an $\IR$-set of $G$ and $\Gamma(G) = \IR(G)$, this implies that $B$ is the unique $\Gamma$-set of $G$.~\QED

\medskip
We show next that Conjecture~\ref{conj1} holds for connected claw-free graphs.

\begin{theorem}
\label{t:clawfree2}
If $G$ is a connected, claw-free graph of order~$n \ge 2$, then the following holds. \\ [-26pt]
\begin{enumerate}
\item $\Psi_g^+(G) \ge \frac{1}{2}n$. \1
\item If $G \ne P_3$, then $\Psi_g^-(G) \ge \frac{1}{2}n$.
\end{enumerate}
\end{theorem}
\proof Let $G$ be a connected, claw-free graph of order~$n \ge 2$. Suppose the Maximizer-start enclaveless game is played on $G$. Let $S$ denote the set of all vertices played when the game ends. By definition of the game, the set $S$ is a maximal enclaveless set in $G$. If $\IR(G) \le \frac{1}{2}n$, then analogously as in the proof of Theorem~\ref{t:bound_IR1} we have $|S| = \Psi_g^+(G) \ge \psi(G) \ge n - \IR(G) \ge \frac{1}{2}n$. Hence, we may assume that $\IR(G) > \frac{1}{2}n$, for otherwise the desired result follows. By Theorem~\ref{t:bound_IR2}, $\IR(G) = \frac{1}{2}(n+1)$ and $G \in \cF$. Let $(A,B)$ be the partition associated with $G$. We show in this case we have $\Psi_g^+(G) > \psi(G)$.

By Observation~\ref{ob:relate}, $\Gamma(G) + \psi(G) = n$. Moreover, the complement of every $\Gamma$-set of $G$ is a maximal enclaveless set, and the complement of every $\psi$-set of $G$ is a minimal dominating set. By Corollary~\ref{c:cor1}, the set $B$ is the unique $\Gamma$-set of $G$. These observations imply that the complement of the set $B$, namely the set $A$, is the unique $\psi$-set of $G$. Thus every maximal enclaveless set of $G$ of cardinality~$\psi(G)$ is precisely the set $A$.

We now return to the Maximizer-start enclaveless game played on $G$. If Maximizer plays as his first move any vertex from the set $B$ and thereafter both players play optimally, then the resulting set $S^*$ of moves played during the course of the game contain a vertex of $B$ and is therefore different from the set $A$. Since the set $A$ is the unique $\psi$-set of $G$, this implies that $|S^*| > \psi(G)$. We therefore have that the following inequality chain, where the first inequality, namely $\Psi_g^+(G) \ge |S^*|$, is due to the fact that the first move of Maximizer from the set $B$ may not be an optimal move.
\[
\Psi_g^+(G) \ge |S^*| \ge  \psi(G) + 1 = (n - \Gamma(G)) + 1 = n - \frac{1}{2}(n+1) + 1 = \frac{1}{2}(n+1).
\]

This shows that $\Psi_g^+(G) \ge \frac{1}{2}n$, as desired. Suppose next that $G \ne P_3$ and the Minimizer-start enclaveless game is played on $G$. Let $S$ denote the set of all vertices played when the game ends. By definition of the game, the set $S$ is a maximal enclaveless set in $G$. If $\IR(G) \le \frac{1}{2}n$, then analogously as before we have $|S| = \Psi_g^-(G) \ge \psi(G) \ge n - \IR(G) \ge \frac{1}{2}n$. Hence, we may assume that $\IR(G) > \frac{1}{2}n$, for otherwise the desired result follows. By Theorem~\ref{t:bound_IR2}, $\IR(G) = \frac{1}{2}(n+1)$ and $G \in \cF$. Let $(A,B)$ be the partition associated with $G$.

We show in this case we have $\Psi_g^-(G) > \psi(G)$. Since $G \ne P_3$, we note that there are at least two vertices in the set $B$ at distance at least~$3$ apart in $G$. Thus, whatever the first move is played by Minimizer, Maximizer can always respond by playing as his first move a vertex chosen from the set $B$. Thus, analogously as before, the resulting set of played vertices in the game is different from the set $A$. Recall that upon completion of the game the resulting set is a maximal enclaveless set. Therefore, Maximizer has a strategy to finish the game in at least~$\psi(G) + 1$ moves, implying that $\Psi_g^-(G) \ge \frac{1}{2}(n+1)$.~\QED

\medskip
By Theorem~\ref{t:clawfree2}(a), we note that Conjecture~\ref{conj1} holds for connected claw-free graphs. Moreover by Theorem~\ref{t:clawfree2}(b), we note that Conjecture~\ref{conj3} holds for connected claw-free graphs even if we relax the minimum degree two condition and replace it with the requirement that the graph is isolate-free and different from the path $P_3$.

\end{document}